\documentclass[reqno]{amsart}    
\usepackage{pb-diagram,enumerate}    
\usepackage{amssymb}    
    
\theoremstyle{plain}    
\newtheorem{thm}{Theorem}[section]    
\newtheorem{cor}[thm]{Corollary}    
\newtheorem{lem}[thm]{Lemma}    
\newtheorem{prop}[thm]{Proposition}    
\theoremstyle{remark}    
\newtheorem{rem}[thm]{Remark}    
\theoremstyle{definition}    
\newtheorem{defn}[thm]{Definition}    
\newtheorem{example}[thm]{Example}    
    
\def\al{{\alpha}}    
    
\def\de{{\delta}}    
    
\def\om{{\omega}}    
    
\def\la{{\lambda}}

\def\si{{\sigma}}    
\def\Si{{\Sigma}}    
\def\ga{{\gamma}}    
\def\ep{{\varepsilon}}

\def\phi{{\varphi}}    
\def\Ph{{\Phi}}

\DeclareMathAlphabet{\doba}{U}{msb}{m}{n}

\gdef\mN{\doba{N}}

\gdef\mR{\doba{R}}    
\gdef\mS{\doba{S}}

\gdef\Cl{\doba{C}l}

\def\lamin{\lambda_{\rm min}^+}    
\def\laminm{\lambda_{\rm min}^-}    
\def\Vol{{\mathop{\rm Vol}}}

\def\Aut{{\mathop{\rm Aut}}}    
\def\End{{\mathop{\rm End}}}    
\def\Hom{{\mathop{\rm Hom}}}    
\def\Spin{{\mathop{\rm Spin}}}    
\def\SO{{\mathop{\rm SO}}}    
\def\rank{{\mathop{\rm rank}}}    
\def\ker{{\mathop{\rm ker}}}    
    
\def\geucl{{g_{\rm eucl}}}    
\def\Geucl{{G_{\rm eucl}}}    
\def\Deucl{{D_{\rm eucl}}}    
    
\def\Id{{\mathop{\rm Id}}}    
    
\def\eref#1{(\ref{#1})}    
\def\mand{{\qquad\mbox{and}\qquad}}

\let\<\langle    
\let\>\rangle

\long\def\komment#1{}

\begin{document}    
    
\title{Mass endomorphism and  spinorial Yamabe type problems on conformally flat manifolds}    
\author{Bernd Ammann, Emmanuel Humbert and Bertrand Morel}
\date{\today}    
    
\begin{abstract}    
Let  $M$ be a compact manifold equipped with a Riemannian metric $g$ and a    
spin structure $\si$.    
We let $\lamin (M,[g],\si)= \inf_{\tilde{g} \in [g] } \lambda_1^+(\tilde{g})    
Vol(M,\tilde{g})^{1/n} $ where $\lambda_1^+(\tilde{g})$ is the    
smallest positive eigenvalue of the    
Dirac operator $D$ in the metric $\tilde{g}$.    
A previous result stated that $\lamin(M,[g],\si) \leq  \lamin(\mS^n) =\frac{n}{2} \om_n^{{1/n}}$    
where $\om_n$ stands for the volume of the standard    
$n$-sphere. In this paper, we study this problem for conformally flat    
manifolds of dimension $n \geq 2$ such that $D$ is invertible. E.g. we show that strict inequality holds in dimension    
$n\equiv 0,1,2\mod 4$ if a certain endomorphism does not vanish.    
Because of its tight relations to the ADM mass in General Relativity,    
the endomorphism will be called mass endomorphism. We apply the strict    
inequality to spin-conformal spectral theory and show that    
the smallest positive Dirac eigenvalue attains its infimum inside the enlarged    
volume-$1$-conformal class of $g$.    
\end{abstract}    
\thanks{Ammann was supported by the NSF grant DMS-9810361 of MSRI, Berkeley}       
\maketitle   

\tableofcontents    
    
\section{Introduction}    
    
Let $(M,g,\si)$ be a compact spin manifold of dimension $n \geq 2$.    
For a metric $\tilde g$ in the conformal class $[g]$ of $g$,    
let $\lambda_1^+(\tilde{g})$ be the smallest positive    
eigenvalue of the Dirac operator $D$. Similarly, let    
$\lambda_1^-(\tilde{g})$ be the largest negative    
eigenvalue of $D$.  We define    
$$\lamin(M,[g],\si) = \inf_{\tilde{g} \in [g] } \lambda_1^+(\tilde{g})    
\Vol(M,\tilde{g})^{1/n} $$    
and    
$$\laminm(M,[g],\si) =  \inf_{\tilde{g} \in [g] } |\lambda_1^-(\tilde{g})|\;    
\Vol(M,\tilde{g})^{1/n} $$    
    
It was proven in \cite{ammann:03} that    
$$\lamin(M,[g],\si)>0 \quad\mbox{and}\quad \laminm(M,[g],\si)) >0.$$    
Several works have been devoted to the study of this conformal    
invariant. A non-exhaustive list is    
\cite{hijazi:86,hijazi:91,lott:86,baer:92b,ammann:p03a,ammann:habil}.    
In \cite{ammann:03}, we proved the following result:    
\begin{thm} \label{theo.old}    
Let $(M,g,\si)$ be a compact spin manifold of dimension $n \geq 2$.    
Then    
  $$\lamin(M,[g],\si) \leq \lamin(\mS^n) = \frac{n}{2}\, \om_n^{{1 \over n}}\mand \laminm(M,[g],\si)\leq \lamin(\mS^n),$$    
where $\om_n$ stands for the volume of the standard sphere    
$\mS^n$.     
\end{thm}    
\noindent The strict inequalities    
\begin{eqnarray} \label{str_in}    
\lamin(M,[g],\si) <  \frac{n}{2}\, \om_n^{{1 \over    
    n}} \mand \laminm(M,[g],\si) <  \frac{n}{2}\, \om_n^{{1 \over    
    n}} \;,    
\end{eqnarray}    
have several applications. At first, when $n >  2$, 
together with Hijazi inequality, each one of the two inequalities    
\eref{str_in} implies the existence of a solution of the Yamabe problem. 
This    
problem is a famous problem of conformal geometry which has been solved by    
Aubin \cite{aubin:76} and Schoen \cite{schoen:84}. 
 
Another application of the inequalities \eref{str_in} is the    
solution of a conformally invariant PDE    
which can be read as a nonlinear eigenvalue equation for the Dirac    
operator.    
The nonlinearity involves a critical exponent from the point of view of Sobolev    
embeddings and hence, this PDE cannot be solved by standard    
methods. Moreover,    
from this solution, one can construct a generalized metric    
$\tilde{g}$ (see \cite{ammann:habil,ammann:p03a})    
such that $Vol_{\tilde{g}}(M)=1$ and such that   
$\lambda_1^+(\tilde{g}) =  \lamin(M,[g],\si)$    
(resp. $\lambda_1^-(\tilde{g}) = \laminm(M,[g],\si)$).   
In other words, this proves that   
$  \lamin(M,[g],\si)$ (resp. $\laminm(M,[g],\si)$) is attained by the    
generalized metric $\tilde{g}$.    
    
In this paper    
we are interested in obtaining the inequalities (\ref{str_in}) when $M$ is    
a     
conformally flat manifold of dimension $n \geq 2$ such that $D$ is invertible. In    
this goal,    
we introduce in section 2 the notion of    
\emph{mass endomorphism}. This endomorphism    
corresponds to the constant term in the development of the Green function    
for $D$ near the diagonal with respect to a conformal chart.  
In Remark~\ref{rem.two.ten} we will show that   
the pointwise eigenvalues of the mass endomorphism are all real.  
The mass endomorphism   
plays the same role as the constant term of the Green function    
$\gamma(\,\cdot\,,p)$    
of the Yamabe operator.    
In the Yamabe problem,    
the constant term of $\gamma$    
can be interpreted as the mass \cite{arnowitt.deser.misner:62}    
of the asymptotically    
flat manifold $(M\setminus\{p\},\gamma(\,\cdot\,,p)^{4/{(n-2)}}g)$    
(see also \cite{lee.parker:87}).    
This is why we use the name \emph{mass endomorphism}.    
Schoen shows in \cite{schoen:84} that the positivity    
of the mass implies the solution of the Yamabe problem.    
In this paper, the eigenvalues of the mass endomorphism play the same role    
as the mass in Yamabe problem.   
Namely, we obtain the following result:    
    
\begin{thm} \label{main2}    
Let $(M,g,\si)$ be a conformally flat    
compact spin manifold of dimension $n \geq 2$ with $\ker(D)= \{0\}$.     
Assume that the mass endomorphism (see next section) possesses a positive     
(resp.\ negative)    
eigenvalue. Then    
$$ \lamin(M,[g],\si) \; \quad(\hbox{resp.\ } \laminm(M,[g],\si) ) < \lamin(\mS^n) = \frac{n}{2}\, \om_n^{{1 \over    
    n}}.$$     
\end{thm}    
    
\noindent    
Assume that $n \not\equiv 3 \mod 4$,     
then the spectrum of the Dirac operator and   
the pointwise spectrum of the mass endomorphism are symmetric  
(see Subsection~\ref{subsec.symm}).    
In particular,    
  $$\lamin(M,[g],\si)= \laminm(M,[g],\si).$$    
This implies 
\begin{thm} \label{main1}    
Let $(M,g,\si)$ be a conformally flat    
compact spin manifold of dimension $n \geq 2$ with $\ker(D)= \{0\}$.    
Assume that  $n \not\equiv 3 \mod 4$ and that the mass endomorphism $\alpha$  
is not identically zero. Then  
$$\lamin(M,[g],\si) < \lamin(\mS^n) = \frac{n}{2}\, \om_n^{{1 \over    
    n}}\mand \laminm(M,[g],\si)< \lamin(\mS^n) = \frac{n}{2}\, \om_n^{{1 \over  n}}.$$    
\end{thm}
     
\noindent This is no longer true if $n\equiv 3\mod 4$.   
In Example~\ref{ex.rpn}    
we study the real projective spaces $\mR P^{4k+3}$. Here the mass    
endomorphism is a non-vanishing multiple of the identity section,    
hence has constant sign which depends on the spin structure.    
Furthermore, the two spin structures    
$\si_+$ and $\si_-$ on $\mR P^{4k+3}$ satisfy    
  $$\laminm(\mR P^{4k+3},g_0,\si_+)>\lamin(\mR    
  P^{4k+3},g_0,\si_+)= {n\over 2}\left({\om_n\over 2}\right)^{1\over    
  n},$$    
  $$\lamin(\mR P^{4k+3},g_0,\si_-)>\laminm(\mR    
  P^{4k+3},g_0,\si_-)= {n\over 2}\left({\om_n\over 2}\right)^{1\over    
  n}.$$    
Hence, as predicted by Theorem~\ref{main2} either $\lamin$ or $\laminm$   
is smaller than $(n/2)\om_n^{1/n}$. As a remark, if $\min(\lamin(\mR    
  P^{4k+3},g_0,\si_-),\laminm(\mR    
  P^{4k+3},g_0,\si_-) )<
(n/2)\om_n^{1/n}$ then the infimum in the definition  
of $\la_{\rm min}^\pm$ would be attained by a metric of non-constant   
curvature. It is then natural to think that we  cannot obtain the strict
inequality \eref{str_in} for $\lamin$ and $\laminm$  
for all manifolds of dimension $n\equiv 3 \mod 4$. 
    
In order to prove inequalities (\ref{str_in}) for arbitrary    
conformally flat manifolds of dimension $n \geq 2$ such that $D$ is    
invertible, then one has to find some nonzero eigenvalues of the mass    
endomorphism. Some questions arise naturally. At first, on the sphere    
$\mS^n$, the mass endomorphism is null. 
Otherwise, we could apply Theorem~\ref{main2}.     
One may wonder if the mass endormorphism is not always zero. The answer is   
no since  as  mentioned above,  the    
projective spaces give examples for which the mass endomorphism is a nonzero   
multiple of the identity.   
  
Recall once again that the constant term of the Green    
function of the Yamabe operator can be interpreted as the mass of an    
asymptotically flat manifold; according to the positive mass theorem     
it is positive unless is $M$    
is conformally diffeomorphic to $\mS^n$. Then, one could hope to find a     
result of    
the same  type for the mass endomorphism. However,    
we show  that this is false in general. Namely, we show in Section~\ref{sec.inv} that    
the mass endomorphism of flat tori always vanishes.

\section{The mass endomorphism}    
    
In the following, we assume that $(M,g,\si)$ is a conformally flat manifold    
of dimension $n \geq 2$. We also assume that $D$ is invertible, i.e. that     
$\ker(D) = \{0\}$.    
The mass endomorphism is defined as the constant term in the Green    
function     
for~$D$.    
In this section, we give a precise definition and some properties of    
the mass endomorphism.    
    
\subsection{Trivialization of the spinor bundle}    
Let $\rho:U\subset (M,g)\to V\subset (N,g_N)$ be a    
conformal map preserving orientation and spin structure.    
We write    
$g= f^2 \rho^*g_N$. Then according to \cite{hitchin:74,hijazi:86} there is    
a fiberwise isomorphism $$    
\begin{array}{rccl}    
  {\rho_*}:& \Si M|_U &\rightarrow& \Si N|_V \\    
  & \downarrow&& \downarrow\\    
  & U & \rightarrow & V,    
\end{array}    
$$    
such that $D_N \rho_*  \phi =  f  \rho_* D\phi$  for all    
$\phi\in \Gamma(\Sigma M)$ and    
$f^{-\frac{n-1}{2}}  \rho_*$    
is an isometry, where $D$ and $D_N$ denote the    
Dirac operators on $M$ and $N$. The most important case we will use    
is that $N$ is Euclidean space, i.e.\ $\rho$ is a conformal chart of $(M,g)$  
preserving orientation and the spin structure.    
In this case the above map yields a well chosen trivialization of the    
spinor bundle.    
  
The definition of the mass endomorphism will be done by working with a   
conformal  charts.    
For simplicity we will first define it in the special case that $g$ is   
flat in a neighborhood of a given base point, and then extend it   
to the general case.  

\subsection{Green function for the Dirac operator}

Let $\pi_1,\pi_2:M\times M\to M$ be the projection to the first and    
second component. Then we define    
  $$\Si M\boxtimes \Si M^*:=\pi_1^*(\Si M)\otimes (\pi_2^*(\Si M))^*,$$    
i.e. it is the bundle whose fiber over $(x,y)$ is given by    
$\Hom(\Si_y M,\Si_x M)$. Let $\Delta:=\{(p,p)\,|\,p\in M\}$ be the    
diagonal.

\begin{defn}\label{def.Green}    
A smooth section $G_D:M\times M \setminus \Delta\to \Si M\boxtimes \Si M^*$    
that is locally integrable on $M\times M$ 
is called \emph{the Green function} for the Dirac operator $D$    
if in the sense of distributions    
$$D_x(G_D(x,y)) = \delta_y Id_{\Si_yM}$$    
In other words, we have for any $y\in M$, $\psi_0\in \Si_y M$, and     
$\phi\in \Gamma(\Si M)$    
  $$\int \<G_D(x,y)\psi_0,D^*\phi(x)\>\,dx= \<\psi_0,\phi(y) \>.$$     
\end{defn}    
    
For $k:M\times M\setminus \Delta \to  \Si M\boxtimes \Si M^*$ we denote    
the corresponding integral operator by $P_k$, i.e.\ we define for     
$\phi\in \Gamma(\Sigma M)$     
  $$P_k(\phi)(x):=\int k(x,y)\phi(y)\,dy.$$    
If $k$ is smooth on $M\times M\setminus \Delta$ and locally integrable
on $M\times M$, then the integral exists in the principal value sense of 
distributions. $P_k$ uniquely determines $k$.    
In this notation $P_{G_D}$ is the inverse of $D:L_1^2(\Si M)\to L^2(\Si M)$.     
Hence, the Green function is unique.    
    
\begin{rem} \label{green_eucl}      
Analogously, one can define the Green function of Euclidean space as a section     
$\Geucl:\mR^n\times \mR^n\setminus \Delta\to \Si M\boxtimes \Si M^*$    
such that $P_\Geucl$ is the inverse of $\Deucl:L_1^2(\Si M)\to L^2(\Si M)$.    
Again we have unicity.    
One easily checks that     
  $$\Geucl(x,y)=-{1\over \om_{n-1}}\, {x-y\over |x-y|^n}\cdot.$$    
\end{rem}

In our construction of the test spinor,    
we will need the asymptotics of Green functions close to the diagonal.    
This is provided by the following  
proposition.    
    
\begin{prop}\label{expan}\    
Assume that the metric is flat near $y \in M$, and let $\rho$ be an    
\emph{isometric} chart.    
Then, the Green function $G_D$ for~$D$ exists, and     
in the above trivialization       
has the following expansion when $x$ tends    
to $y$    
$$\omega_{n-1} \,G_D(x,y) \psi_0= - \frac{x-y}{|x-y|^n}  \cdot \psi_0 +    
v(x,y)\psi_0$$    
where $D_x v(x,y)(\psi_0)=0$ on a neighborhood of $y$.    
\end{prop}    
\noindent {\bf Proof:}    
Let $G_{\Deucl}$ be the Green function for $\Deucl$ on $\mR^n$ given by    
Remark~\ref{green_eucl}.    
We take a cut-off function $\eta$ with support in $B_y(\de)$    
which is equal to $1$ on $B_y(\de/2)$ where $\de>0$ is a small number.    
We set $\Ph(x)= \eta(x)  G_{\Deucl}(x,y) \cdot \psi_0$ where $\psi_0$ is    
constant. The spinor  $\Ph$ is harmonic on $B_y(\de/2)\setminus  \{y\}$.    
We extend $\Phi$ by zero, and obtain a smooth spinor on $M\setminus  \{y\}$.    
As $D\Phi|_{B_y(\de/2)}\equiv 0$, we see that $D\Phi$ extends to    
a smooth spinor on $M$, denoted by $\Psi$.    
Since $D$ is assumed to be invertible, there exists a smooth spinor field    
$\zeta$ such that     
$D\zeta= -\Psi$.    
We then define     
$G_D(x,y)\psi_0:=\Phi(x)-\zeta(x)$.      
It is easy to see that $G_D$ is the Green function for $D$ and has the    
development described above.    
As $\Psi$ and hence $\zeta$ depend smoothly on $y$,   
it is clear that $G_D(x,y)$  
depends smoothly on $y$ outside of the diagonal.  
\qed

\begin{lem}[Conformal change and Green functions]\label{lem.conf.Green}    
Let $\rho:(M,g_M,\si)\to (N,g_N,\si_N)$ be a conformal diffeomorphism 
preserving   
the orientation and the spin structure  
and write $g_M=f^2\rho^* g_N$. Let $G_M$ resp. $G_N$ be the Green function    
on $M$ resp. $N$, then    
  $$G_N(\rho(x),\rho(y))= f^{n-1}(y)  \rho_{*,x} \circ G_M(x,y) \circ \rho_{*,y}^{-1}.$$    
\end{lem}    
    
\proof{}    
We already know that $D_N=f \rho_* D\rho_*^{-1}$.    
This implies     
$$P_{G_N}= D_N^{-1}=\rho_*D^{-1}\rho_*^{-1}f^{-1}=\rho_*P_{G_M}\rho_*^{-1}f^{-1},$$    
Hence,     
  $$G_N(\rho(x),\rho(y))\,d\mu_N(\rho(y))= f^{-1}(y)  \rho_{*,x} \circ {G_M}(x,y) \circ \rho_{*,y}^{-1}\,d\mu(y).$$    
With $d\mu=f^nd\mu_N$ this implies the lemma.    
\qed    
    
In particular, the previous proposition and the lemma imply that the     
Green function $G_D$ for $D$ exists on any conformally flat manifold for
which $D$ is invertible.

\begin{example}    
If $\phi_{\la}$ is an eigenspinor for the    
eigenvalue $\la$ of $D$, and if $G_D$ is the unique     
Green function for $D$, then    
$$ \int_M G_D(x,y) \phi_{\la} (y)\,dy =     
\frac{1}{\la} \phi_{\la}(x) $$    
\end{example}    
  
\begin{example}\label{example.green.RS}    
Let $\rho:S^n \setminus \{N\}\to \mR^n$ be the stereographic projection.    
Then     
  $$g_{S^n}= \left({2\over |\rho(x)|^2+1}\right)^2\,\rho^*\geucl.$$     
Let $G_{S^n}$ be the Green function on $S^n$.     
An obvious modification of Lemma~\ref{lem.conf.Green} tells us    
that    
  $$\rho_* G_{S^n}(x,y)\rho_*^{-1}=      
    \left({2\over |\rho(y)|^2+1}\right)^{1-n} \Geucl(\rho(x),\rho(y)).$$    
\end{example}    
  
\begin{example}\label{example.green.SS}    
Let $\rho:S^n \to S^n$ be a Moebius transformation of $S^n$  
with $g_{S^n}= f^2\rho^*(g_{S^n})$. Then  
  $$G_{S^n}(x,y)= f^{n-1}(y)  \rho_{*,x} \circ G_{S^n}(x,y) \circ \rho_{*,y}^{-1}.$$    
\end{example}

\noindent In the sequel  the following self-adjointness result for $G_D$ will be important.    
    
\begin{prop} \label{GDsa}    
Let $V$ be an open set of $M$ in which $g$ is flat. Then, for all $x \not=    
y \in V$, $G_D(x,y)^*=G_D(y,x)$. In other words,    
$$\<G_D(x,y) \psi_y, \phi_x\>=\<\psi_y, G_D(y,x) \phi_x\>$$    
for all $\psi_y \in \Si_y M$ and $\phi_x \in \Si_x M$.    
\end{prop}    
    
\noindent {\bf Proof:} 
We have   
  $$P_{G_D}^*=\left(D^{-1}\right)^*=\left(D^*\right)^{-1}= D^{-1}= P_{G_D}.$$  
As the operator uniquely determines the kernel, this implies the proposition.   
\qed  
  
\begin{prop}\label{prop.conf.masse}  
Assume that we have two metrics   
$g_1$ and $g_2$ on $M$, $g_1=f^2 g_2$. We assume that both metrics    
are flat in a neighborhood $U$ of $y$. We define $v_1$ and $v_2$ as above.    
Then   
  $$v_{1}(x,x)=f^{1-n} v_2(x,x)$$    
\end{prop}    
\proof{}   
Let $G^1$ (resp. $G^2$) be the Green function for $g_1$ (resp. $g_2$). We identify spinors via $\rho_*$. By    
Lemma \ref{lem.conf.Green}, $G^1$ and $G^2$ are related by the following    
formula: for all $x,y \in M$,    
\begin{equation}\label{eq.G.one.two}    
  G^1(x,y)=  f^{1-n}(y) G^2(x,y)    
\end{equation}    
Let $h_1:(U,g_1)\to \mR^n$ and $h_2:(U,g_2)\to \mR^n$ be isometric embeddings.    
According to Liouville's theorem  
the conformal diffeomorphism $h_2\circ h_1^{-1}: h_1(U)\to h_2(U)$  
extends to a M\"obius  transformation of the Alexandrov compactification $S^n$ of $\mR^n$.    
Because of Examples~\ref{example.green.RS} and \ref{example.green.SS}   
this implies that  
  $$G_{\rm eucl}^1(x,y)=  f^{1-n}(y) G_{\rm eucl}^2(x,y).$$    
Subtracting this from \eref{eq.G.one.two}, and taking the limit $x\to y$    
one obtains the desired formula.    
\qed

\subsection{Definition and first properties of the mass endomorphism}    
    
\noindent The mass endomorphism is defined as the constant term of the     
Green function for $D$ with respect to a conformal chart.     
Let us make this precise.    
\begin{defn} Let $(M,g)$ be a compact manifold which is conformally flat     
on a neighborhood of $y\in M$. Choose a metric $\tilde g\in [g]$     
that is flat on a neighborhood of $y$ and such that $\tilde g_y=g_y$.    
Let $G_D$ be the Green function for $D$. Then we define     
the \emph{mass endomorphism} as    
\[ \alpha_y : \left|  \begin{array}{ccc}    
\Si_y(M) & \rightarrow & \Si_y(M) \\    
\psi_0 & \mapsto & v(y,y)(\psi_0)    
\end{array}    
\right. \]    
where $v$ is as in the previous paragraph with respect to  $\tilde g$.    
\end{defn}    
    
Because of Proposition~\ref{prop.conf.masse} this definition does not depend     
on the choice of $\tilde g$.

\begin{prop} \label{mass_prop}    
For each $y \in M$, the mass endomorphism $\al_y$ is linear and self-adjoint.  
\end{prop}    
    
\noindent {\bf Proof:}    
Let $y \in M$ and  $\psi_0, \phi_0 \in \Si_y(M)$. We have    
$$\al_y(\psi_0) = \lim_{x \to y} \om_{n-1} G_D(x,y) \psi_0 +    
\frac{x-y}{|x-y|^n}  \cdot \psi_0$$    
It follows immediately that $\al_y$ is linear.    
Taking the limit when $x \to y$     
one gets that    
$$\<\al_y(\psi_0), \phi_0\>= \<\psi_0, \al_y(\phi_0)\>$$    
\qed    
    
\begin{rem}\label{rem.two.ten}    
Proposition~\ref{mass_prop} immediately implies that the mass endomorphism    
has only real eigenvalues.     
\end{rem}    
    
\begin{rem}\label{rem.quat.comm}    
If $\dim M\equiv 2,3,4\,\mod 8$, then the spinor bundle     
carries a quaternionic structure, i.e.\ a basepoint-preserving, parallel,     
complex anti-linear map     
$Q:\Si M\to \Si M$ with $Q^{-1}$ commuting with the Clifford multiplication.    
As a consequence $Q$ commutes with $D$, $G_D$ and with $\alpha$.    
\end{rem}    
    
\subsection{Examples}    
\begin{example}[Flat tori]    
Let $(M,g)$ be an $n$-dimensional flat torus.   
It carries $2^n$ spin structures.    
For one spin structure, the so-called \emph{trivial} spin structure,   
we have a $\rank(\Sigma M)=2^{[n/2]}$-dimensional    
space of parallel sections. All other spin structures admit no    
non-trivial parallel spinors. Because $M$ is scalar-flat, the kernel of $D$    
consists exactly of the parallel spinors, in particular $D$ is invertible  
for all non-trivial spin structures $\chi$.    
Because translations  act spin-isometric on $(M,g,\chi)$,   
the Green function $G_D$ satisfies    
$G_D(x,y)=G_D(x-y,0)$. Also $D_x(G_D(y,x))=-\delta(x-y)$, hence    
$G_D(x,y)=-G_D(y,x)$. Therefore, all terms of even order in    
the development of $G_D$ have to vanish. In particular, the mass    
endomorphism vanishes.    
\end{example}    
    
\begin{example}[Real Projective Spaces]\label{ex.rpn}    
Besides $\mR P^1=S^1$, the only real projective spaces that are orientable and spin are    
$\mR P^{4k+3}$ with $k\in \mN$.    
The space $M=\mR P^{4k+3}$ carries exactly    
two spin structures. The universal covering $\pi:S^{4k+3}\to \mR P^{4k+3}$    
induces a push-forward of the spinor bundles,    
which is a fiberwise isomorphism    
$\pi_*:\Si_pS^{4k+3}\to \Si_{\pi(p)}\mR P^{4k+3}$. One calculates    
\begin{equation}\label{eq.G.rpn}
G_D^{\mR P^{4k+3}}(\pi{x},\pi{y})\circ \pi_*= \pi_*G_D^{S^{4k+3}}(x,y)+    
    \pi_*G_D^{S^{4k+3}}(x,-y),
\end{equation}    
where $-y$ denotes the antipodal point of $y$.   
Stereographic projection based in $y\not\in \{x,-x\}$    
defines a conformal chart containing $x$ and $-x$.   
Example~\ref{example.green.RS} implies that   
$G_D^{S^{4k+3}}(x,-x)\neq 0$. Hence, the mass endomorphism    
$\al_{\pi(x)}$ of $\mR P^{4k+3}$ does not vanish anywhere on $\mR P^{4k+3}$.    
  
The group of orientation preserving   
isometries fixing $\pi(x)$ is $\SO(4k+3)$. After passing to the double    
cover $\Spin(4k+3)$, we obtain a $\Spin(4k+3)$-action on $\Si S^{4k+3}$, that  
pushes down to a $\Spin(4k+3)$-action on $\Si S^{4k+3}$ which    
commutes with the Dirac operator. Hence, this action also commutes with the   
mass endomorphism, and as the $\Spin(4k+3)$-action on    
$\Si_{\pi(x)}\mR P^{4k+3}$ is irreducible, the mass endomorphism is a    
constant multiple of the identity \cite[Prop.~I.5.15]{lawson.michelsohn:89},
\cite[section~1.5]{friedrich:00}.
If one changes the spin structure, then the second summand in \eref{eq.G.rpn}
changes its sign.  
Hence, the sign of the mass endomorphism   
depends on the choice of the spin structure, which are denoted by $\si_+$ and  
$\si_-$.    
\end{example}    
    
\subsection{Endomorphisms generating symmetries}\label{subsec.symm}    
    
The aim of this section is to show    
that if $n \not\equiv 3 \mod 4$, then there is an automorphism    
$\Gamma(\Aut_{\mR}(\Sigma M))$ that anticommutes with the Dirac operator.   
This result is well-known, see for example \cite[1.7]{friedrich:00},    
\cite[Prop. 5]{dahl:p03}.   
As a consequence, one sees that it also anticommutes with the    
mass endomorphism.    
   
Let $(W,\gamma)$ be an irreducible    
complex representation of the Clifford algebra    
of an Euclidean vector space $V$ of dimension $n$. After fixing an    
orientation on $V$, one can define    
  $$\omega_\gamma:= \gamma(e_1\cdot \ldots \cdot e_n)\in\End(W)$$    
where $e_1,\ldots,e_n$ is an oriented orthonormal basis on $V$, and where    
$\gamma$ denotes Clifford multiplication. One easily calculates    
  $$\omega^2:=(-1)^{n(n+1)/2}.$$    
As a consequence the eigenvalues of $\omega$ are contained in $\{-1,1\}$    
if $n\equiv 0,3\mod 4$, and they are contained in $\{-i,i\}$    
if $n\equiv 1,2\mod 4$.    
    
\def\makecase#1{{\ \newline\noindent \it #1}\\[2mm]}    
\makecase{The case $n\equiv 0\mod 2$.}    
If $n$ is even, then Clifford multiplication by a vector $v\in V$    
anticommutes with $\omega$. Hence, if $\alpha$ is an eigenvalue $\om$,    
then so is $-\alpha$.    
One immediately obtains the well-known lemma.    
\begin{lem}    
If $n$ is even, then    
$\omega$ is a complex-linear automorphism of $W$ anticommuting    
with Clifford-multiplication.    
\end{lem}    
Indeed, it can even be shown that up to automorphisms there is only one    
irreducible representation in even dimensions $n$ and that any endomorphism    
anticommuting with Clifford multiplication with vectors is a multiple of    
$\omega$.    
    
\makecase{The case $n\equiv 1\mod 4$.}    
The question arises, whether there is a similar endomorphism    
if $n$ is odd.    
In this case, Clifford multiplication with a vector commutes with    
$\omega$. Hence, by Schur's lemma $\omega$ has only one eigenvalue.    
For $n\equiv 1\mod 4$ we have either $\omega=i\,\Id$ or $\omega=-i\,\Id$, and it    
can be shown, that there is exactly one irreducible    
representation of $V$ with $\omega=i\,\Id$ denoted by $(W^i,\gamma^{i})$,    
and one with $\omega=-i\,\Id$, denoted by $(W^{-i},\gamma^{-i})$.    
If we replace the complex structure on $W^{\pm i}$, by its complex    
conjugate one, then this is again    
a representation of the (real) Clifford algebra of $V$.    
Obviously, $\om$ changes sign by conjugation.    
Hence, there is a conjugate linear isomorphism of Clifford representation    
  $\alpha:W^{i} \to W^{-i}.$    
Another way  to modify the structure of $(W^i,\gamma^i)$ is    
to reverse the sign of Clifford multiplication by vectors.    
Namely, we define a Clifford multiplication $\tau^i:\Cl(V)\to \End(W)$    
as $\tau^{i}(X):=-\gamma^{i}(X)$  
for all vectors $X$ in $V$.    
Again, we calculate that the sign of $\omega$ changes if we replace    
$\tau^{i}$ by $\gamma^i$, and there is a    
complex linear isomorphism of vector spaces    
  $\beta:W^{-i} \to W^{i}$    
with $\beta\circ \ga^{-i}(X)=  \tau^{i}(X)\circ \beta$ for any vector $X\in V$.    
Hence, $\nu:=\beta\circ \alpha:W^i\to W^i$ is a conjugate linear automorphism    
of vector spaces, and for vector $X$ we have    
  $$\gamma^i(X)\circ \nu = -\tau^{i}(X)\circ \beta\circ \alpha =    
    - \beta\circ\gamma^{i}(X)\circ\alpha = - \nu\circ    
    \gamma^{i}(X).$$    
A similar endomorphism of $W^{-i}$ is given by $\alpha\circ \beta$.    
We have proven the following lemma.    
\begin{lem}    
If $n\equiv 1\mod 4$, then    
there is a real vector space automorphism of $W^{\pm i}$ anticommuting with    
Clifford multiplication by vectors. The automorphism is conjugate linear.    
\end{lem}    
    
\makecase{The case $n\equiv 3\mod 4$.}    
The case $n\equiv 3\mod 4$ is different. Again,    
we have $\omega= \Id$ or $\omega= -\Id$,    
and there is exactly one irreducible    
representation $(W^\pm,\gamma^\pm)$ in each case.    
However, conjugation does not exchange the representations, and    
the sign of $\omega$ is    
invariant under real automorphisms.    
Hence, an isomorphism as in the above lemma    
cannot exist.    
    
\makecase{For all $n\not\equiv 3\mod 4$.}    
Any automorphism $\nu$ anticommuting with Clifford multiplication, commutes    
with bivectors $X\cdot Y$ where $X,Y\in V$. As the Lie algebra of $\Spin(n)$    
is generated by elements of that form, such an isomorphisms $\nu$ is    
$\Spin(n)$-equivariant. We obtain the well-known    
\begin{prop}    
If $n\equiv 0,1,2\mod 4$, then there is a real vector bundle isomorphism    
$\nu:\Si M\to \Si M$    
anticommuting with Clifford multiplication by vectors,    
complex linear if $n$ is even, and conjugate linear if $n\equiv 1\mod 4$.    
Furthermore, $\nu$ is parallel.    
\end{prop}    
    
It follows that it anti-commutes with the Dirac operator, the Green function    
and the mass endomorphism.    
    
\begin{cor}[Well-known, e.g.\ \cite{atiyah.patodi.singer:75}]    
The spectrum of the Dirac operator is symmetric in dimension $n\not \equiv 3\mod 4$.    
\end{cor}    
    
\begin{cor}    
The pointwise spectrum of the mass endomorphism is symmetric in    
dimension $n\not \equiv 3\mod 4$, i.e.\ if $\lambda$ is an eigenvalue of the    
mass endomorphism $\alpha_x$ for an $x\in M$,    
then $-\lambda$ is also an eigenvalue of $\al_x$ with the same    
multiplicity.    
\end{cor}

\begin{cor}    
If $\dim M=2$, then the mass endomorphism $\alpha$ vanishes.    
\end{cor}    
\begin{proof}    
The spectrum of $\al_y$ is symmetric and real. As $\alpha$ commutes    
with the quaternionic multiplication $Q$, the eigenspaces of $\alpha$ are    
quaternionic vector spaces. $\dim \Si_y M=2$ implies that $\al_y=r_y\Id$    
for $r_y\in\mR$. As the spectrum of $\al_y$ is symmetric, we obtain    
$r_y=0$.    
\end{proof}

\section{The estimates}\label{sec.est}    
    
\noindent Let $\psi\in\Gamma(\Sigma M)$ and define    
    
$$J(\psi)=\frac{\Big(\int_M|D\psi|^{\frac{2n}{n+1}}v_g\Big)^\frac{n+1}{n}}{\int_M\Re    
  e\<D\psi,\psi\>v_g}$$    
The first named author proved in \cite{ammann:03} that    
\begin{eqnarray}\label{deflamin1}    
\lamin(M,[g],\si)=\inf_\psi J(\psi)    
\end{eqnarray}    
where the infimum is taken over the set of smooth spinor fields for which    
$$\left(\int_M    
\Re e  \<D\psi,\psi\>v_g \right)>0\;.$$    
By adjusting some signs appropriately, one obtains by the same reasoning that    
\begin{eqnarray}\label{deflamin2}    
\laminm(M,[g],\si)=\inf_\psi - J(\psi)    
\end{eqnarray}    
where the infimum is taken over the set of smooth spinor fields for which    
$$\left(\int_M    
\Re e  \<D\psi,\psi\>v_g \right)<0\;.$$    
These two facts will be helpful to prove the following theorem    
    
\begin{thm}\label{maintheorem}    
Assume that there exists on $(M,g)$ a conformal metric $g_p\in [g]$ which is flat in a neighborhood of a point $p$.    
If there exists on $M\setminus \{p\}$ a spinor field $\psi$ satisfying    
\begin{itemize}    
  \item $D\psi=0$    
  \item $\psi$ admits the following development near the point $p$: $$\psi=\frac{x}{r^n}    
\cdot \psi_0 +\psi_1+\theta\;$$ where $\psi_0$ and $\psi_1$ are two spinors    
of $\Sigma_p M$ such that    
$$\Re e(\<\psi_0, \psi_1\>) <0\qquad(\text{resp. } \Re e(\<\psi_0,    
\psi_1\>)>0)\;,$$    
and where  $\theta=O(r)$ is an harmonic spinor field smoothly defined in a    
neighborhood of $p$.    
\end{itemize}    
Then we have    
$$ \lamin(M,[g],\si) <  \frac{n}{2}\, \om_n^{{1 \over    
    n}} \qquad \Big(\text{resp. } \laminm(M,[g],\si) <  \frac{n}{2}\, \om_n^{{1 \over    
    n}}\Big)\;.$$    
\end{thm}    
    
\begin{proof} The proof is based on a suitable choice of a test spinor field $\psi_\varepsilon$ to estimate    
$J(\psi_\varepsilon)$. We let    
    
$$f(r):=\frac{1}{1+r^2}\;.$$    
    
\noindent As a first step, we consider for    
a given $\Phi$ the spinor defined by    
  $$\varphi^{\pm}(x):= {f(|x|)}^{{n \over 2}}(1 \mp x)\cdot \Phi \qquad \forall x \in \mR^n.$$    
One may compute that    
  $$D(\varphi^{\pm}) = \pm n f(|x|) \varphi^{\pm}\;.$$    
Since the Euclidean space $(\mR^n \setminus  \{ 0\} ,4 f^2 g_{eucl})$ and the standard sphere    
$(\mathbb{S}^n,g_0)$ are isometric and using the    
conformal covariance of $D$, it is well known  that  there exists a    
natural map    
\[ m \left| \begin{array}{ccc}    
\Gamma (\Si(\mR^n \setminus  \{ 0\})) & \to& \Gamma(\Si(\mathbb{S}^n)) \\    
\phi & \mapsto & m(\phi)    
\end{array} \right. \]    
such that for all vector field $\phi \in \Gamma (\Si(\mR^n \setminus  \{    
0\}))$, we have    
  $$ m(D(\phi)) = 2 f D_{\mathbb{S}^n} (m(\phi))$$    
As one can check the spinor field $m(\varphi^{+})$ is    
a Killing spinor on $\mathbb{S}^n$ to the Killing constant $-{1 \over 2}$,   
 whereas    
$m(\varphi^{-})$ is a Killing spinor on $\mathbb{S}^n$ to the Killing    
constant $+{1 \over 2}$.    
One gets    
that    
$$J_{\mR^n}(\varphi^\pm)= \pm \frac{n}{2}  \om_n^{{1 \over    
    n}}$$    
where $J_{\mR^n}$ is the functional $J$ written on $\mR^n$.    
Now, fix  $\varepsilon>0$ and let $(x_1,\cdots,x_n)$ be local coordinates    
on a neighborhood $U$ of $p$ in $M$. On $U$, we trivialize the spinor    
bundle via parallel transport. Using this trivialization, one may define    
for $x \in M$,    
    
$$\varphi^\pm_\varepsilon := \eta_0 \varphi^\pm(\frac{x}{\ep})$$    
where $\eta_0$ is a cut-off function equal to $1$ on $B(p,\delta)$    
($\delta$ is a small number). 

We would like to have test spinors $\psi_\pm$ with    
$\pm\int \<D\psi_\pm,\psi_\pm\>>0$ for which the \emph{strict} inequalities    
  $$\pm J(\psi_\pm) < \frac{n}{2}\, \om_n^{{1 \over n}}$$    
hold.    
For a given $\ep>0$ we set    
  $$\xi:=\varepsilon^{\frac{1}{n+1}}\qquad \ep_0:=\frac{\xi^n}{\varepsilon} f(\frac{\xi}{\varepsilon})^{\frac{n}{2}}\;.$$    
The test spinor we use here is the    
following:    
    
\begin{equation}\label{defpsieps}    
\psi_\ep^\pm:=    
\begin{cases}    
f(\frac{r}{\varepsilon})^{\frac{n}{2}}(1\mp\frac{x}{\varepsilon})\cdot\psi_0\mp    
\varepsilon_0 \psi_1& \text{ if } r\leq\xi\;, \\[5mm]    
\mp\varepsilon_0(\psi-\eta \theta)+\eta f(\frac{\xi}{\varepsilon})^{\frac{n}{2}}\psi_0 & \text{ if } \xi\leq\ r\leq 2\xi\;,\\[5mm]    
\mp \varepsilon_0 \,\psi&\text{ if } r\geq 2\xi\;,    
\end{cases}    
\end{equation}    
where $r=|x|$, where $\eta$ is a cut-off function which equals    
to $1$ on $B(p,\xi)$, which is zero on the complement    
of $B(p,2\xi)$ and which statisfies    
  $$|\nabla\eta|\leq\frac{2}{\xi}.$$    
    
\noindent Note that $\psi^\pm_{\ep}$ is continuous on $M$.    
\begin{rem}    
It should be pointed out that this choice of $\xi$ is arbitrary in the    
following sense.  The proof of Theorem    
\ref{maintheorem} still holds for any choice of $\xi=\varepsilon^q$ for    
$q\in ]\frac{n-1}{n(n+1)}\, ,\,\frac{1}{n}[$.   
\end{rem}    
We can assume without loss of generality that $|\psi_0|=1$. Since $\psi$    
and $\theta$ are harmonic near $p$, we have    
    
\begin{equation}\label{Dpsieps}    
D\psi_\varepsilon^\pm=    
\begin{cases}    
\pm\frac{n}{\varepsilon}f(\frac{r}{\varepsilon})^{\frac{n}{2}+1}(1\mp\frac{x}{\varepsilon})\cdot\psi_0&    
\text{ if } r\leq\xi\;,\\[5mm]    
\pm\varepsilon_0\nabla\eta\cdot\theta+f(\frac{\xi}{\varepsilon})^{\frac{n}{2}}\nabla\eta\cdot\psi_0 & \text{ if } \xi\leq\ r\leq 2\xi\;,\\[5mm]    
0&\text{ if } r\geq 2\xi\;.    
\end{cases}    
\end{equation}    
    
Therefore, since $|(1\mp\frac{x}{\varepsilon})\cdot\psi_0|^2=(1+\frac{r^2}{\varepsilon^2})|\psi_0|^2=f(\frac{r}{\varepsilon})^{-1}$, we have    
    
\begin{equation}\label{normDpsieps}    
|D\psi_\varepsilon^\pm|^{\frac{2n}{n+1}}=    
\begin{cases}    
\left[\frac{n}{\varepsilon}f(\frac{r}{\varepsilon})^{\frac{n+1}{2}}\right]^{\frac{2n}{n+1}}=n^{\frac{2n}{n+1}}\varepsilon^{-\frac{2n}{n+1}}f(\frac{r}{\varepsilon})^n& \text{ if } r\leq\xi\;,\\[5mm]    
|\pm\varepsilon_0\nabla\eta\cdot\theta+f(\frac{\xi}{\varepsilon})^{\frac{n}{2}}\nabla\eta\cdot\psi_0|^{\frac{2n}{n+1}} & \text{ if } \xi\leq\ r\leq 2\xi\;,\\[5mm]    
0&\text{ if } r\geq 2\xi\;.    
\end{cases}    
\end{equation}    
    
In the following, the notation $C$ will stand for positive    
constants (eventually depending on the dimension $n$ but not on $\ep$)    
which can differ from    
one line to another. Equation (\ref{normDpsieps}) yields the following estimates:    
$$\int_{B(p,\xi)}|D\psi_\varepsilon^\pm|^{\frac{2n}{n+1}}=    
\varepsilon^{n-\frac{2n}{n+1}}n^{\frac{2n}{n+1}}\int_{B(p,    
  \frac{\xi}{\ep})}    
f^n    
\leq \varepsilon^{n-\frac{2n}{n+1}}n^{\frac{2n}{n+1}}\int_{\mathbb{R}^n}f^n$$    
    
and    
    
\begin{eqnarray*}    
\int_{B(p,2\xi)\setminus B(p,\xi)}|D\psi_\varepsilon^\pm|^{\frac{2n}{n+1}}&\leq&C\,\int_{B(p,2\xi)\setminus B(p,\xi)} |\,\varepsilon_0\nabla\eta\cdot\theta|^{\frac{2n}{n+1}}+C'\,\int_{B(p,2\xi)\setminus B(p,\xi)}|f(\frac{\xi}{\varepsilon})^{\frac{n}{2}}\nabla\eta\cdot\psi_0|^{\frac{2n}{n+1}}\\[3mm]    
&\leq&C\,\varepsilon^{\frac{n(2n-1)}{n+1}}+    
C'\,\varepsilon^{\frac{n(2n-1)}{n+1}} \leq C\,\varepsilon^{\frac{n(2n-1)}{n+1}}\;,    
\end{eqnarray*}    
    
since $\varepsilon_0\leq C\,\varepsilon^{n-1}$,    
$|\nabla\eta|\leq\frac{2}{\xi}$, $Vol(B(p,2\xi)\setminus B(p,\xi)) \leq C \xi^n$    
and $|\theta|\leq C\,\xi$ on    
$B(p,2\xi)$, as well as $f(\frac{\xi}{\varepsilon})^{\frac{n}{2}}\leq C\,\varepsilon^{\frac{n^2}{n+1}}$.    
    
Therefore    
\begin{equation}\label{estnormdpsi}    
\left(\int_M|D\psi_\varepsilon^\pm|^{\frac{2n}{n+1}}\right)^\frac{n+1}{n}\leq    
\varepsilon^{n-1}n^2I^{1+\frac{1}{n}}\left[1+C\,\varepsilon^{\frac{n^2}{n+1}}    
\right]= \varepsilon^{n-1}n^2I^{1+\frac{1}{n}}\left[1+o(\ep^{n-1})    
\right]\;,    
\end{equation}    
    
where $I:=\int_{\mathbb{R}^n}f^n$.    
    
\noindent If we set    
$$\nu:= \<\psi_0,\psi_1\>\;,$$ we have    
\begin{eqnarray*}    
\Re e\<D\psi_\varepsilon^\pm, \psi_\varepsilon^\pm\>_{\vert B(p,\xi)}&=&\Re    
e\<\pm\frac{n}{\varepsilon}f(\frac{r}{\varepsilon})^{\frac{n}{2}+1}(1\mp\frac{x}{\varepsilon})\cdot\psi_0\,,f(\frac{r}{\varepsilon})^{\frac{n}{2}}(1\mp\frac{x}{\varepsilon})\cdot\psi_0\mp\varepsilon_0    
\, \psi_1\>\\    
&=&\pm\frac{n}{\varepsilon}f(\frac{r}{\varepsilon})^n-\frac{n}{\varepsilon}\,\varepsilon_0\,    
f(\frac{r}{\varepsilon})^{\frac{n}{2}+1}\,\Re    
e(\nu) \pm\frac{n}{\varepsilon}\varepsilon_0\,    
f(\frac{r}{\varepsilon})^{\frac{n}{2}+1}\,\Re    
e\<\frac{x}{\varepsilon}\cdot\psi_0,\psi_1\>\,    
\end{eqnarray*}    
and hence, since by symmetry the last term vanishes when    
integrating over $B(p,\xi)$, we have    
\begin{eqnarray*}    
\int_{B(p,\xi)}\Re e\<D\psi_\varepsilon^\pm, \psi_\varepsilon^\pm\>    
&=&n\,\varepsilon^{n-1}\Big[\pm\int_{B(p,\frac{\xi}{\varepsilon})}f(r)^n-\Re e(\nu) \varepsilon_0\,\int_{B(p,\frac{\xi}{\varepsilon})}    
f(r)^{\frac{n}{2}+1}\Big]\;.    
\end{eqnarray*}    
Moreover,    
$$\int_{B(p,\frac{\xi}{\varepsilon})}f(r)^n=I-\omega_{n-1}\int_{\frac{\xi}{\varepsilon}}^{+\infty}r^{n-1}    
f(r)^n\mathrm{d}r\;,$$    
Since    
$$\int_{\frac{\xi}{\varepsilon}}^{+\infty}r^{n-1}    
f(r)^n\mathrm{d}r \leq \int_{\frac{\xi}{\varepsilon}}^{+\infty} r^{-(n+1)}    
\mathrm{d}r \leq C \,\varepsilon^{\frac{n^2}{n+1}}\;$$    
and since $\ep_0\sim \ep^{n-1}$    
when $\ep \to 0$, we have, for $\Re e(\nu)<0$,    
\begin{equation}    
  \label{fff}    
 \int_{B(p,\xi)}\Re e\<D\psi_\varepsilon^+, \psi_\varepsilon^+\>\geq    
n\,\varepsilon^{n-1}\Big[I-C_0\,\Re e(\nu)    
\varepsilon^{n-1}+o(\ep^{n-1})\Big]\;,    
\end{equation}    
and for $\Re e(\nu)>0$,    
\begin{equation}    
  \label{vvv}    
-\int_{B(p,\xi)}\Re e\<D\psi_\varepsilon^-, \psi_\varepsilon^-\>\geq    
n\,\varepsilon^{n-1}\Big[I+C_0\,\Re e(\nu)    
\varepsilon^{n-1}+o(\ep^{n-1})\Big]\;,    
\end{equation}    
where    
$$C_0= \int_{\mR^n}    
f(r)^{\frac{n}{2}+1}$$    
We also have    
\begin{eqnarray*}    
\Re e\<D\psi_\varepsilon^\pm, \psi_\varepsilon^\pm\>_{\vert B(p,2\xi)\setminus B(p,\xi)}&=&\Re e\<\pm\varepsilon_0\nabla\eta\cdot\theta+f(\frac{\xi}{\varepsilon})^{\frac{n}{2}}\nabla\eta\cdot\psi_0 ,\mp\varepsilon_0(\psi-\eta \theta)+\eta f(\frac{\xi}{\varepsilon})^{\frac{n}{2}}\psi_0\>\\    
&=&\Re e\<\pm\varepsilon_0\nabla\eta\cdot\theta+f(\frac{\xi}{\varepsilon})^{\frac{n}{2}}\nabla\eta\cdot\psi_0 ,\mp\varepsilon_0\,\psi\>    
\end{eqnarray*}    
since $\Re e\<\nabla\eta\cdot\theta,\theta\>=0\;$, $\Re e\<\nabla\eta\cdot\psi_0,\psi_0\>=0$    
and    
$$\Re e\<\nabla\eta\cdot\psi_0,\theta\>+\Re e\<\nabla\eta\cdot\theta,\psi_0\>=\Re e\<\nabla\eta\cdot\psi_0,\theta\>-\Re e\overline{\<\nabla\eta\cdot\psi_0,\theta\>}=0\;.$$    
Now, we write that on $ B(p,2\xi) \setminus  B(p,\xi)$, $\ep_0 |\psi| \leq    
C \ep^{n-1} \xi^{-(n-1)}$, $\ep_0 | \nabla \eta| |\theta| \leq C \ep^{n-1}$ and    
$f(\frac{\xi}{\varepsilon})^{\frac{n}{2}}|\nabla\eta| \leq C    
\ep^{n-1}$. This leads to    
    
\begin{eqnarray*}    
\Re e \<D\psi_\varepsilon^\pm, \psi_\varepsilon^\pm\>_{\vert    
  B(p,2\xi)\setminus  B(p,\xi)}& \leq & C \ep^{2n-2} \xi^{1-n} \;,    
\end{eqnarray*}    
which yields    
\begin{eqnarray}\label{lastterm}    
\int_{B(p,2\xi)\setminus B(p,\xi)}\Re e\<D\psi_\varepsilon^\pm,    
\psi_\varepsilon^\pm\>=O(\varepsilon^{2n-2+\frac{1}{n+1}})= o(\ep^{2n-2})\;.    
\end{eqnarray}    
Therefore, since    
\begin{eqnarray*}    
\Re e\<D\psi_\varepsilon^\pm, \psi_\varepsilon^\pm\>_{\vert M\setminus B(p,2\xi)}&=&0\;,    
\end{eqnarray*}    
from Equations \eqref{fff}, \eqref{vvv}, \eqref{lastterm},    
we have for $\Re e(\nu)<0$,    
$$\int_M\Re e\<D\psi_\varepsilon^+, \psi_\varepsilon^+\>\geq    
n\,\varepsilon^{n-1}I\Big[1-C_0\,\Re e(\nu)    
\varepsilon^{n-1}+o(\ep^{n-1})\Big]\;,$$    
and for $\Re e(\nu)>0$,    
$$-\int_M\Re e\<D\psi_\varepsilon^-, \psi_\varepsilon^-\>\geq    
n\,\varepsilon^{n-1}I\Big[1+C_0\,\Re e(\nu)    
\varepsilon^{n-1}+o(\ep^{n-1})\Big]\;,$$    
    
\noindent Together with \eqref{estnormdpsi}, we then have, for  $\Re e(\nu)<0$,    
$$J(\psi_\varepsilon^+)\leq n\,I^{\frac{1}{n}}\,\frac{1+o(\ep^{n-1})}    
{1-C_0\,\Re e(\nu)    
\varepsilon^{n-1}+o(\ep^{n-1})} \;,$$    
and for  $\Re e(\nu)>0$,    
$$-J(\psi_\varepsilon^-)\leq    
n\,I^{\frac{1}{n}}\,\frac{1+o(\ep^{n-1})}    
 {1+C_0\,\Re e(\nu) \varepsilon^{n-1}+ o(\ep^{n-1})}    
\;.$$    
For $\ep$ small enough, we obtain    
$$J(\psi_\varepsilon^+) < n I^{\frac{1}{n}} \hbox{ and } -J(\psi_\varepsilon^-)    
 < n I^{\frac{1}{n}}$$    
Recall now the following fact: let $p$ be any point of the round sphere    
$\mathbb{S}^n$. Then $\mathbb{S}^n\setminus  \{p\}$ is isometric to    
$\mathbb{R}^n$ with the metric $$g_\mathbb{S}=4f^2g_\mathrm{eucl}\;.$$    
    
Therefore $$2^nI=\int_{\mathbb{R}^n}2^nf^n=\omega_n\;,$$    
which yields, for $\Re e(\nu)<0$ (resp. for $\Re e(\nu)>0$)    
$$J(\psi_\varepsilon^+)<\frac{n}{2}\,\omega_n^{\frac{1}{n}}\qquad    
(\text{resp.}\quad    
-J(\psi_\varepsilon^-)<\frac{n}{2}\,\omega_n^{\frac{1}{n}})\;.$$    
Hence, by \eqref{deflamin1} and \eqref{deflamin2}, the proof of Theorem    
\ref{maintheorem} is now complete.    
\end{proof}

\section{Proofs of Theorems~\ref{main2} and \ref{main1}}\label{sec.inv}    
    
Let $p \in M$. Up to a conformal change of metric, we may assume that $g$    
is flat near $p$. Assume that the mass endomorphism $\al_p$ possesses a    
non-zero eigenvalue $\nu$. Let $\psi_0 \in \Si_p(M) $ be    
an eigenvector associated to $\nu$. Then, we set    
$$\psi=  -\om_{n-1} G_D(x,p)    
\psi_0$$    
The spinor field $\psi$ then satisfies the assumptions of    
Theorem \ref{maintheorem} with $\psi_1= \nu \psi_0$.    
Theorem \ref{maintheorem} implies that if $\nu>0$ then $\lamin(M,g,\theta) < \frac{n}{2}\, \om_n^{{1 \over    
    n}}$ and if $\nu<0$, $\laminm(M,g,\theta) < \frac{n}{2}\, \om_n^{{1    
    \over  n}}$. This proves Theorem~\ref{main2}. Now, if $n \not\equiv 3    
\mod 4$, the    
spectrum of the mass endomorphism is symmetric and hence if $\al_p \not=    
0$, $\nu$ can be chosen positive or negative.    
This proves Theorem~\ref{main1}.

\vspace{1cm}    
Author's addresses:    
\nopagebreak    
\vspace{5mm}\\    
\parskip0ex    
{\obeylines    
Bernd Ammann, Emmanuel Humbert, and Bertrand Morel,    
Institut \'Elie Cartan BP 239    
Universit\'e de Nancy 1    
54506 Vandoeuvre-l\`es -Nancy Cedex    
France    
}    
\vspace{0.5cm}    
    
E-Mail:    
{\tt ammann@iecn.u-nancy.fr, humbert@iecn.u-nancy.fr and morel@iecn.u-nancy.fr}    
    
\end{document}